\documentclass[a4paper,12pt]{amsart}
\usepackage{amssymb}
\usepackage{dsfont}
\usepackage{ifthen}
\usepackage{graphicx}
\usepackage{color}

\newtheorem{thm}{Theorem}
\newtheorem{cor}{Corollary}
\newtheorem{lem}{Lemma}

\newtheorem{que}{Question}
\newtheorem{rem}{Remark}
\newtheorem{prob}{Problem}
\newtheorem{conj}{Conjecture}
\theoremstyle{definition}
\newtheorem{example}[equation]{Example}

\newcounter {own}
\def\theown {\thesection       .\arabic{own}}

\newenvironment{pf}[1][]{%
 \vskip 3mm
 \noindent
 \ifthenelse{\equal{#1}{}}%
  {{\slshape Proof. }}%
  {{\slshape #1.} }%
 }%
{\qed\bigskip}

\newcounter{alphabet}
\newcounter{tmp}
\newenvironment{Thm}[1][]{\refstepcounter{alphabet}%
\bigskip%
\noindent%
{\bf Theorem \Alph{alphabet}}%
\ifthenelse{\equal{#1}{}}{}{ (#1)}%
{\bf .} \itshape}{\vskip 8pt}

\makeatletter
\newcommand{\Ref}[1]{\@ifundefined{r@#1}{}{\setcounter{tmp}{\ref{#1}}\Alph{tmp}}}
\makeatother

\newcounter{minutes}
\divide\time by 60
\newcounter{hours}
\multiply\time by 60
\addtocounter{minutes}{-\time}

\newcommand{\IC}{{\mathbb C}}
\newcommand{\ID}{{\mathbb D}}

\def\be{\begin{equation}}
\def\ee{\end{equation}}

\newcommand{\bee}{\begin{enumerate}}
\newcommand{\eee}{\end{enumerate}}

\newcommand{\blem}{\begin{lem}}
\newcommand{\elem}{\end{lem}}
\newcommand{\bthm}{\begin{thm}}
\newcommand{\ethm}{\end{thm}}
\newcommand{\bcor}{\begin{cor}}
\newcommand{\ecor}{\end{cor}}
\newcommand{\beg}{\begin{example}}
\newcommand{\eeg}{\end{example}}
\newcommand{\bques}{\begin{que}}
\newcommand{\eques}{\end{que}}
\newcommand{\begs}{\begin{examples}}
\newcommand{\eegs}{\end{examples}}
\newcommand{\bdefe}{\begin{defin}}
\newcommand{\edefe}{\end{defin}}
\newcommand{\bprob}{\begin{prob}}
\newcommand{\eprob}{\end{prob}}
\newcommand{\bei}{\begin{itemize}}
\newcommand{\eei}{\end{itemize}}

\newcommand{\bcon}{\begin{conj}}
\newcommand{\econ}{\end{conj}}
\newcommand{\bcons}{\begin{conjs}}
\newcommand{\econs}{\end{conjs}}
\newcommand{\bprop}{\begin{propo}}
\newcommand{\eprop}{\end{propo}}
\newcommand{\br}{\begin{rem}}
\newcommand{\er}{\end{rem}}
\newcommand{\brs}{\begin{rems}}
\newcommand{\ers}{\end{rems}}
\newcommand{\bo}{\begin{obser}}
\newcommand{\eo}{\end{obser}}
\newcommand{\bos}{\begin{obsers}}
\newcommand{\eos}{\end{obsers}}
\newcommand{\bpf}{\begin{pf}}
\newcommand{\epf}{\end{pf}}
\newcommand{\ba}{\begin{array}}
\newcommand{\ea}{\end{array}}
\newcommand{\beq}{\begin{eqnarray}}
\newcommand{\beqq}{\begin{eqnarray*}}
\newcommand{\eeq}{\end{eqnarray}}
\newcommand{\eeqq}{\end{eqnarray*}}

\newcommand{\ra}{\rightarrow}

\newcommand{\ds}{\displaystyle}

\def\cc{\setcounter{equation}{0}   
\setcounter{figure}{0}\setcounter{table}{0}}

\begin{document}

\bibliographystyle{amsplain}

%

\title[Univalent functions defined by differential inequalities]{On classes of meromorphic locally univalent functions defined by differential
inequalities}

\thanks{
File:~\jobname .tex,
          printed: \number\day-\number\month-\number\year,
          \thehours.\ifnum\theminutes<10{0}\fi\theminutes}

\author[S. K. Lee]{See Keong Lee}
\address{S. K. Lee,
School of Mathematical Sciences,
Universiti Sains Malaysia,
11800 USM Penang, Malaysia.}

\email{sklee@usm.my}

\author[S. Ponnusamy]{Saminathan Ponnusamy
}

\address{S. Ponnusamy, Department of Mathematics,
Indian Institute of Technology Madras,
Chennai-600 036, India. }

\email{samy@iitm.ac.in}

\author[K.-J. Wirths]{Karl-Joachim Wirths}
\address{K.-J. Wirths, Institut f\"ur Analysis und Algebra, TU Braunschweig,
38106 Braunschweig, Germany.}
\email{kjwirths@tu-bs.de}


\begin{abstract}
In this article we consider functions meromorphic in the unit disk. We give an elementary proof
for a condition that is sufficient for the univalence of such functions which also contains some
known results. We include few open problems for further research.
\end{abstract}


\maketitle
\pagestyle{myheadings}
\markboth{S. K. Lee, S. Ponnusamy and K.-J. Wirths}{Univalent functions defined by differential inequalities}
\cc

\section{Preliminaries and Main Results}
We denote the unit disk by $\ID=\{z\in \IC:\, |z| < 1\}$ and let
$${\mathcal A}=\{f:\,\mbox{ $f$ is analytic in $\ID$, $f(0)=f'(0)-1=0$}\}.
$$
The family ${\mathcal S}$ of univalent functions in $\mathcal A$ together with its many subfamilies, for which the image domains
have special geometric properties, have been investigated in details. See \cite{Duren:univ, Go}. Throughout, $\mathcal{B}$ denotes
the class of functions $\omega$, analytic in $\ID$ such that $|\omega(z)|\leq 1$ for $z\in \ID$.
The well known inequality $ \sum_{k=0}^{\infty}|c_k|^2\leq 1$ for $\Omega(z)= \sum_{k=0}^{\infty}c_kz^k \in \mathcal{B}$ will be used to get the
proof of Theorem \ref{th4}. Recently,  in \cite[Theorem 2(b)]{PW-2017}, the second and the third authors proved among other things
the following result which extends the earlier known result for analytic functions.

\begin{Thm}\label{Th1}
Let $f$ be meromorphic in $\ID$ such that $f(0)=f'(0)-1=0$. If for all $z\in\ID$  the inequality
\begin{equation}\label{f5}
\left|\frac{z}{f(z)} - z\left(\frac{z}{f(z)}\right)'- 1\right| < \lambda, ~z\in\ID,
\end{equation}
is valid for some $\lambda \in (0,1]$, then $f$  is univalent in $\ID$.
\end{Thm}

The proof of Theorem \Ref{Th1} in \cite{PW-2017} was elegant and was also different from the other known methods.
In this article, we shall consider slightly more general situation. For $0<\lambda$ and $\mu \in\IC$ such that $|1-\mu|<\lambda$, we consider the family ${\mathcal U}(\lambda,\mu)$ of meromorphic functions $f$ satisfying the inequality
$$  |U_f(z)-\mu| <\lambda ~\mbox{ in $\ID$},
$$
where
\be\label{OS-eq2}
U_f(z):=\left (\frac{z}{f(z)} \right )^{2}f'(z)=\frac{z}{f(z)} -z\left (\frac{z}{f(z)} \right )', \quad z\in\ID.
\ee
Note that the center $1$ has been replaced by $\mu$. First we consider the problem of determining conditions on $\lambda$ and $\mu$ so that
functions in  ${\mathcal U}(\lambda,\mu)$ are univalent in $\ID$.

As with the case of analytic functions, for notational simplicity, we let ${\mathcal U}(\lambda):={\mathcal U}(\lambda,1)$, ${\mathcal U}:={\mathcal U}(1)$. In the analytic case, it was well-known that ${\mathcal U}\subsetneq {\mathcal S}$ (see \cite{Aks58,AksAvh70}).

\bthm\label{th1}
Let $\lambda \in (0,1]$ and $|1-\mu|<\lambda$. All members of the  family $\mathcal{U}(\lambda,\mu)$ are functions meromorphic and  locally univalent in $\ID$ if and only if $|\mu| \geq \lambda.$
\ethm
\bpf
It is a simple exercise to see that if $f\in {\mathcal U}(\lambda,\mu)$, then $f(z)\neq 0$ in $\ID \backslash\{0\}$, because otherwise
$f(z_0)=0$ for some $z_0\in\ID \backslash\{0\}$ which would then imply that
$$f(z)=c_m(z-z_0)^m+\cdots  \quad (c_m\neq 0,~m\geq 1)
$$
in a neighborhood of $z_0$ so that
$$U_f(z) 
=\frac{mz_0}{c_m(z-z_0)^{m+1}} +\cdots
$$
implying that $U_f(z)$ has a pole of order $m+1$ which would clearly be a contradiction to the fact that $U_f(z)$ is bounded. Moreover,
for $|\mu|\geq \lambda$, the inequality $|U_f(z)-\mu|<\lambda$ for $z\in \ID$ implies $|U_f(z)|>0$ for $z\in \ID$ and hence, in either way $f'(z)\neq 0$ for
$z\in \ID$ and $f(z)\neq 0$ in $\ID \backslash\{0\}$.

To prove the other direction of our assertion, we let $f(z) = z + \sum_{n=2}^{\infty} a_n z^n$ in ${\mathcal U}(\lambda,\mu)$ and consider
$$ \lambda \Omega (z)=
\frac{z}{f(z)} -z\left (\frac{z}{f(z)} \right )'-\mu=1-\mu +(a_3-a_2^2) z^2 +\cdots =1-\mu + O(z^2)
$$
as $z\ra 0$, where   $\Omega$ is analytic in $\ID$ such that $\Omega(0)=(1-\mu)/\lambda$,
$\Omega'(0)=0$ and $|\Omega (z)|<1$ in $\ID$. For simplicity, we let $\Omega(0)=a$. Then
\[
\phi(z)=\frac{\Omega(z)-a}{1-\overline{a}\Omega(z)}\]
is a Schwarz function, i.\,e., $\phi(0)=0$ and $|\phi(z)| < 1$ for $z\in\ID$.
Hence,
\[
\Omega(z)=\frac{\phi(z)+a}{1+\overline{a}\phi(z)}.
\]
Since $\Omega'(0)=(1-|a|^2)\phi'(0)=0$, we have $\phi'(0)=0$. Therefore, we can write $\phi(z)=z^2\omega(z)$,
where $\omega \in \mathcal{B}.$ Consequently,
\be\label{Nf2}
\frac{z}{f(z)} -z\left (\frac{z}{f(z)} \right )'-\mu =\lambda \left (\frac{a+z^2\omega(z)}{1+ \overline{a} z^2\omega(z)}\right ), \quad z\in\ID,
\ee
where $a=(1-\mu)/\lambda$ and $\omega \in \mathcal{B}$. Note that, in the neighborhood of $z=0$, we have the representation
\beqq
\frac{z}{f(z)}&=&\frac{1}{1+ \sum_{n=2}^{\infty} a_n z^{n-1}}\\
& =&1-\left ( \sum_{n=2}^{\infty} a_n z^{n-1} \right ) + \left ( \sum_{n=2}^{\infty} a_n z^{n-1} \right )^2
+\cdots \\
&=& 1-a_2z +(a_2^2 - a_3) z^2 +\cdots .
\eeqq
As with the standard procedure, the integration of the differential equation (\ref{Nf2}) delivers
that each $f\in {\mathcal U}(\lambda,\mu)$ has the representation
\be\label{eq-e1}
f(z)= \frac{z}{1 +cz -\lambda z\int_0^z\frac{(1-|a|^2)\omega(t)}{1+\overline{a}t^2\omega(t)}\,dt},
\ee
where $c=-a_2$, $a=\frac{1-\mu}{\lambda}$ and $\omega \in \mathcal{B}$. Now, we let $\omega(z)=1$, $z\in \ID$, in this
representation and obtain that
$$f(z)= \frac{z}{1 +cz -\lambda z\int_0^z\frac{(1-|a|^2)}{1+\overline{a}t^2}\,dt}.
$$
It follows that there exists a $z_1\in \ID$ such that $f'(z_1)=0$ if and only if
\[
\frac{z_1^2}{1+\overline{a}z_1^2} = \frac{-1}{\lambda(1-|a|^2)}.\]
This is equivalent to
\[ |\lambda(1-|a|^2) + \overline{a}|>1,
\]
or equivalently $|\mu|< \lambda$, which is a contradiction to the local univalency of $f$. Hence, the rest of the assertion is proved.
\epf

\bthm\label{th2}
Let $f$ be meromorphic in $\ID$ such that $f(0)=f'(0)-1=0$. Then $f \in U(\lambda,\mu)$ is univalent in $\ID$ if either {\rm (a)} $|1-\mu|< \lambda \leq 1/2$
or {\rm (b)} $\lambda \in (1/2,1]$ and $|1-\mu|\leq 1- \lambda$.
\ethm
\bpf
By using the representation \eqref{eq-e1}, we can write $f=1/g_{\omega}$ for $f\in \mathcal{U}(\lambda,\mu)$, where
$$g_\omega (z)= \frac{1}{z}\left [1 +cz -\lambda z\int_0^z\frac{(1-|a|^2)\omega(t)}{1+\overline{a}t^2\omega(t)}\,dt \right ].
$$
We see that we have to prove
\be\label{th2-eq1}
0 \neq \frac{g_{\omega}(z_1)-g_{\omega}(z_2)}{z_1 -z_2}= \frac{-1}{z_1z_2} -\frac{\lambda}{z_1-z_2}\int_{z_2}^{z_1}\frac{(1-|a|^2)\omega(t)}{1+\overline{a}t^2\omega(t)}\,dt,
\ee
where  $z_1,z_2 \in \ID\setminus\{0\}$, $z_1\neq z_2$, and $\omega\in \mathcal{B}$.
Since $a=(1-\mu)/\lambda$ and
\[
\left|\frac{(1-|a|^2)\omega(z)}{1+\overline{a}z^2\omega(z)}\right| \leq \frac{1-|a|^2}{1-|a|} = 1+|a|,\quad z\in \ID,
\]
we get
\[\left|\lambda \frac{z_1z_2}{z_1-z_2}\int_{z_2}^{z_1}\frac{(1-|a|^2)\omega(t)}{1+\overline{a}t^2\omega(t)}\,dt\right|
< \lambda (1+|a|) = \lambda + |1-\mu|
\]
and thus, \eqref{th2-eq1} holds whenever  $|1-\mu|\leq 1 - \lambda$.
Hence, $f$ is univalent if $|1-\mu|\leq 1 - \lambda$. Therefore, every $f\in \mathcal{U}(\lambda,\mu)$ is univalent in $\ID$
whenever $|1-\mu|<\lambda \leq \frac{1}{2}$ or $|1-\mu|\leq 1-\lambda$ with $\lambda \in (\frac{1}{2},1]$.
This completes the proof.
\epf

\br
{\rm By using Theorem \ref{th1}, we see that at least for nonnegative real numbers $\mu$ the assertion of Theorem \ref{th2}(b) is  best possible. This follows from the fact that for $\mu \in (1-\lambda,\lambda)$, the family $\mathcal{U}(\lambda,\mu)$ contains a function that is not locally univalent in $\ID$. In order to present a couple of precise functions, we consider the function $f_0$ defined by (see also Problem \ref{prob1})
$$ \frac{z}{f_0(z)} =   1-z +\frac{\lambda (1-a^2)z}{2\sqrt{a}} \log
\left ( \frac{1-(\sqrt{a}/(1+\sqrt{a}))(1-z)}{1+(\sqrt{a}/(1-\sqrt{a}))(1-z)} \right ), \quad a=\frac{1-\mu}{\lambda}.
$$
%
}
\er

\bprob
Do there exist families  $\mathcal{U}(\lambda,\mu)$ consisting of univalent functions besides those mentioned in Theorem {\rm \ref{th2}}?
\eprob

In the following we use the equation
\begin{equation}\label{m2}
\frac{z}{f(z)}-z\left(\frac{z}{f(z)}\right)' - \mu = \lambda \Omega(z),
\end{equation}
where
\begin{equation}\label{m3}
\Omega(z)= \frac{1-\mu}{\lambda} + \sum_{k=2}^{\infty}c_kz^k,
\end{equation}
with $\Omega \in \mathcal{B}$ to get sharp estimates for the coefficients of the representation
\begin{equation}\label{m3-a}
\frac{z}{f(z)}=1+\sum_{k=1}^{\infty}b_kz^k=1-a_2z +(a_2^2 - a_3) z^2 +\cdots.
\end{equation}

\bthm\label{th3}
For $f \in \mathcal{U}(\lambda,\mu)$ of the form \eqref{m3-a} and $k\geq 2$, the inequalities
\[
|b_k|\leq \frac{\lambda}{k-1}\left(1 - \frac{|1-\mu|^2}{\lambda^2}\right)
\]
are valid. These inequalities are best possible.
\ethm
\bpf From \eqref{m2} and \eqref{m3} we derive the identities
\[
b_k(1-k) = \lambda c_k,~k\geq 2.
\]
The well known inequalities $|c_k| \leq  1-|c_0|^2$
for $\Omega(z)= \sum_{k=0}^{\infty}c_kz^k$ with $\Omega \in \mathcal{B}$,  $c_0=(1-\mu)/\lambda$ and  $c_1=0$,
imply the validity of our assertion.

For the proof of the sharpness, we set $\omega (z)=-z^{k-2}$ in \eqref{eq-e1} (with $c_0=a$) and  consider the following
functions $f_k$ for $k\geq 2$:
\be\label{m3-b}
\frac{z}{f_k(z)}=1+cz +\lambda z\int_0^z\frac{(1-|a|^2)t^{k-2}}{1-\overline{a}t^k}\,dt
 = 1+cz +\lambda (1-|a|^2)\sum_{j=2}^{\infty}(\overline{a})^j\frac{z^{(j+1)k}}{(j+1)k-1} .
\ee
Obviously, $f_k \in \mathcal{U}(\lambda,\mu)$, and further we get that the $k$-th coefficient $b_k$
of the function $z/f_k(z)$ satisfies
\[
b_k = \frac{\lambda}{k-1}\left(1 - \frac{|1-\mu|^2}{\lambda^2}\right).
\]
The proof is complete.
\epf

\br
{\rm
In particular, in the case $k=2$ of Theorem {\rm \ref{th3}}, we get the sharp inequality
\[
|b_2|=\left|a_2^2 - a_3\right|\leq \lambda\left(1\,-\,\frac{|1-\mu|^2}{\lambda^2}\right)
\]
for $f(z)=z+\sum_{n=2}^{\infty}a_nz^n \in \mathcal{U}(\lambda,\mu).$
}
\er


The well known inequality $ \sum_{k=0}^{\infty}|c_k|^2\leq 1$ for $\Omega(z)= \sum_{k=0}^{\infty}c_kz^k \in \mathcal{B}$ will be used to get the
proof of Theorem \ref{th4}.

\bthm\label{th4}
 If $f\in \mathcal{U}(\lambda,\mu)$ is of the form \eqref{m3-a},
 then
\[
\sum_{k=2}^{\infty}|b_k|^2(k-1)^2\leq \lambda^2\left(1 - \frac{|1-\mu|^2}{\lambda^2}\right).
\]
This inequality is best possible.
\ethm
\bpf
In view of the relations \eqref{m2} and \eqref{m3}, the assumption gives
$$\left |1-\mu -\sum_{k=2}^{\infty}(k-1)b_k z^k\right |\leq \lambda, \quad z\in\ID,
$$
from which the desired inequality follows because $\sum_{k=0}^{\infty}|c_k|^2\leq 1$ for $\Omega(z)= \sum_{k=0}^{\infty}c_kz^k \in \mathcal{B}$. Thus,
it remains to prove the assertion of the sharpness.
To that end, we consider the functions $f_k, k\geq 2$, given by \eqref{m3-b}. Then their Taylor expansions are given by
\[
\frac{z}{f_k(z)}=1+cz +\lambda (1-|a|^2)\sum_{j=0}^{\infty}\overline{a}^j\frac{z^{(j+1)k}}{(j+1)k-1}.
\]
Hence, in these cases we get
\[
\sum_{k=2}^{\infty}|b_k|^2(k-1)^2=\lambda^2(1-|a|^2)^2\sum_{j=0}^{\infty}|a|^{2j}\,= \lambda^2\left(1 - \frac{|1-\mu|^2}{\lambda^2}\right).\]
This completes the proof of the sharpness.
\epf

Finding sharp estimates for the Taylor coefficients of the functions in $U(\lambda, \mu)$ turned out to be a challenge.
As a first result in this direction we prove the next result which extends \cite[Theorem 4]{PW-2017}.

\bthm\label{th5}
If $f \in \mathcal{U}(\lambda,\mu)$ is analytic in the disk $\ID_{p}=\{z: \,|z|<p\}, \, p\in (0,1],$ and $a=(1-\mu)/\lambda$,  then the inequality
\[
|a_2|\leq  A_2:=\frac{1}{p}+\lambda \int_0^p\frac{(1-|a|^2)dt}{1-|a|t^2} = \left \{\ba{rl}
\ds \frac{1}{p}+\frac{\lambda (1-|a|^2)}{2\sqrt{|a|}}\log \left ( \frac{1+p\sqrt{|a|}}{1-p\sqrt{|a|}} \right ) & \mbox{for $a\neq 0$}\\
\ds \frac{1}{p}+\lambda p & \mbox{for $a=0$,}
\ea
\right .
\]
is valid. This estimate is best possible for $\mu \in (1-\lambda, 1].$
\ethm
\bpf We assume on the contrary that $|a_2|>A_2$. In other words, we can assume that there exists an
$r\in (0,1)$ such that
\[
|a_2| = \frac{A_2}{r}=\frac{1}{pr}\left ( 1 +\lambda p\int_0^p\frac{(1-|a|^2)dt}{1-|a|t^2}\right ).
\]
Using Brouwer's fixed point theorem, we shall prove that  then the function
\[
f(z)= \frac{z}{1 - a_2z -\lambda z\int_0^z\frac{(1-|a|^2)\omega(t)}{1+\overline{a}t^2\omega(t)}\,dt},
\]
has a pole in the disk $\overline{\ID}_{rp}:=\{z:\,|z|\,\leq\,pr\}$. To that end, we consider the function
\[
F(z)=\frac{1}{a_2}\left (1 - \lambda z\int_0^z\frac{(1-|a|^2)\omega(t)}{1+\overline{a}t^2\omega(t)}\,dt\right )
\]
and we show that it has a fixed point in the disk $\ID_{pr}$.
For $|z|\leq pr$, we get
\[
|F(z)| \leq \frac{\left(1 + \lambda pr\int_0^{pr}\frac{(1-|a|^2)\,dt}{1-|a|t^2}\right)pr}{1 + \lambda p\int_0^p\frac{(1-|a|^2)\,dt}{1-|a|t^2}}\,< pr.
\]
Since $F$ is a continuous function that maps the convex compact set $\overline{\ID_{pr}}$ into itself, Brouwer's fixed point theorem implies that $F$ has a fixed point in $\overline{\ID_{pr}}$ which is a contradiction to the initial assumptions of Theorem \ref{th5}.
Hence $|a_2|\leq  A_2$ is valid.

Concerning the sharpness, we see that for the numbers $\mu$ in question, the quantity $(1-\mu)/\lambda$ is nonnegative. We choose $\omega(z)=-1$ in the representation formula for $f\in \mathcal{U}(\lambda,\mu)$, and we get that the function $f=f_0$, where
\be\label{LPQ-eq2}
f_0(z)=\frac{z}{1-z\left(\frac{1}{p}+\lambda \int_0^p\frac{(1-a^2)dt}{1-at^2}\right)+\lambda z\int_0^z\frac{(1-a^2)dt}{1-at^2}},\,\,a=\frac{1-\mu}{\lambda},
\ee
which is analytic in $\ID_{p}$ and achieves equality in the estimate of our theorem.
\epf

%
%
%

Note that the function $f_0$ given by \eqref{LPQ-eq2}  takes the form
\beqq
\frac{z}{f_0(z)}
 &=& \left \{\ba{ll}
\ds 1-z \left (\frac{1}{p} +\frac{\lambda (1-a^2)}{2\sqrt{a}}\log \left ( \frac{1+p\sqrt{a}}{1-p\sqrt{a}} \right ) \right) \\
~~~~~~~~~~~~~~~ \ \ \ \ \ds +\,   \frac{\lambda (1-a^2)z}{2\sqrt{a}}\log \left ( \frac{1+\sqrt{a}z}{1-\sqrt{a}z} \right )
  & \mbox{for $a\neq 0$}\\
\ds 1-((1/p)+\lambda ) z +\lambda z^2 & \mbox{for $a=0$,}
\ea
\right .
\eeqq
which may be simplified as
\beqq
\frac{z}{f_0(z)}
&=& \left \{\ba{ll}
\ds 1-\frac{z}{p} +\frac{\lambda (1-a^2)z}{2\sqrt{a}} \log \left ( \frac{1+\sqrt{a}z}{1-\sqrt{a}z} ~\frac{1-p\sqrt{a}}{1+p\sqrt{a}} \right )
  & \mbox{for $a\neq 0$}\\
\ds 1-((1/p)+\lambda ) z +\lambda z^2 & \mbox{for $a=0$,}
\ea
\right .\\
&=& \left \{\ba{ll}
\ds 1-z +\frac{\lambda (1-a^2)z}{2\sqrt{a}} \log
\left ( \frac{1-(\sqrt{a}/(1+p\sqrt{a}))(p-z)}{1+(\sqrt{a}/(1-p\sqrt{a}))(p-z)} \right )
  & \mbox{for $a\neq 0$}\\
\ds 1-((1/p)+\lambda ) z +\lambda z^2 & \mbox{for $a=0$.}
\ea
\right .
\eeqq

In the case of $p=1$,   we then ask in particular the following.

\bprob\label{prob1}
Suppose that $f \in \mathcal{U}(\lambda,\mu)$ is analytic in the unit disk $\ID$, and $a=(1-\mu)/\lambda >0$. Is
$$ \frac{z}{f(z)} \prec \frac{z}{f_0(z)} =   1-z +\frac{\lambda (1-a^2)z}{2\sqrt{a}} \log
\left ( \frac{1-(\sqrt{a}/(1+\sqrt{a}))(1-z)}{1+(\sqrt{a}/(1-\sqrt{a}))(1-z)} \right )?
%
$$
\eprob

Note that this problem has been solved in \cite{OPW} for the case of  $a=0$, i.e. for $\mu =1$ .
\br
{\rm
If  $\mu \notin (1-\lambda, 1]$, it is possible to get an implicit  sharp upper estimate for $|a_2|$ of
functions $f \in \mathcal{U}(\lambda,\mu)$ analytic in the disk $\ID_{p}$ in the following way:
We consider again
\[F(z)=\frac{1}{a_2}\left ( 1 - \lambda z\int_0^z\frac{(1-|a|^2)\omega(t)}{1+\overline{a}t^2\omega(t)}\,dt\right )
\]
and assume that there exists a number $r<1$ such that
\[ |a_2| = \frac{1}{pr}\max\left\{\left|1 - \lambda z\int_0^z\frac{(1-|a|^2)\omega(t)}{1+\overline{a}t^2\omega(t)}\,dt\right| :\,|z|=p,~\omega\in\mathcal{B}\right\}.
\]
Then we use the continuity of the function $F$ on the disk $\overline{\ID_{pr}}$, the inclusion $F(\overline{\ID_{pr}})\subset \overline{\ID_{pr}}$, and Brouwer's fixed point theorem to see that $F$ has a fixed point in $\overline{\ID_{pr}}$. This contradicts the assumption that $f$ is analytic in the disk $\ID_{p}$. Hence,
\[|a_2|\,\leq\,\frac{1}{p}\max\left\{\left|1 - \lambda z\int_0^z\frac{(1-|a|^2)\omega(t)}{1+\overline{a}t^2\omega(t)}\,dt\right| : \,|z|=p,~\omega\in\mathcal{B}\right\}.
\]
To prove the sharpness of this inequality, we choose $z_0$, $|z_0|=p$, and $\omega_0\in \mathcal{B}$ such that

\vspace{8pt}

$\ds
\left|1 - \lambda z_0\int_0^{z_0}\frac{(1-|a|^2)\omega_0(t)}{1+\overline{a}t^2\omega_0(t)}\,dt\right|
=
$
\[
\max\left\{\left|1 -\lambda z\int_0^z\frac{(1-|a|^2)\omega(t)}{1+\overline{a}t^2\omega(t)}\,dt\right| :\, |z|=p, \,\omega\in\mathcal{B}\right\}.
\]
Now, let
\[
a_2= \frac{1}{z_0} - \lambda \int_0^{z_0}\frac{(1-|a|^2)\omega_0(t)}{1+\overline{a}t^2\omega_0(t)}\,dt,
\]
and consider
\[
f_0(z)= \frac{z}{1 - a_2z - \lambda z\int_0^z\frac{(1-|a|^2)\omega_0(t)}{1+\overline{a}t^2\omega_0(t)}\,dt}.
\]
The function $f_0$ shows that the above estimate is sharp.
}
\er

\bprob
Calculate the above maximum.
\eprob

%
%
%
%
\br
{\rm
Let $f(z)=z/g_{\omega}(z)\in \mathcal{U}(\lambda,\mu)$. See \eqref{eq-e1}. Since for any $\omega \in \mathcal{B}$ there exists a positive constant $C$ such that
\[
|g_{\omega}(z_1)- g_{\omega}(z_2)| \leq  C |z_1 - z_2|, ~ z_1,z_2 \in \ID,
\]
the function $g_{\omega}$ is uniformly continuous in $\ID$. Therefore, it has a continuous extension to $\overline{\ID}$. This fact implies that $f$ has a continuous extension to $\overline{\ID}\setminus\{z:\, g_{\omega}(z)=0\}$. Hence, it makes sense to ask for the univalence of this continuous extension. From the proof of Theorem \ref{th2}, it is obvious that this extension is univalent if $|1-\mu| <\lambda \, \leq 1/2$ or if for $\lambda \in (1/2, 1]$ the strict inequality $|1-\mu| < 1 -\lambda$ is valid.

On the other hand, Theorems \ref{th1} and \ref{th2}  have the consequence that for $\lambda \in (1/2, 1]$ the classes $\mathcal{U}(\lambda,\lambda)$ contain the interesting univalent slit mappings
\be\label{LPQ-eq1}
f(z)\,=\,\frac{z}{1 - z\left(\frac{1}{p}\,+\,\lambda \int_0^p\frac{(1-a^2)dt}{1-at^2}\right)\,+\,\lambda z\int_0^z\frac{(1-a^2)dt}{1-at^2}},~a=\frac{1-\lambda}{\lambda}.
\ee
These functions have a pole at $z=p$, and their derivatives vanish at $z=1$ and $z=-1$.  We conjecture that possibly these classes
and those functions deserve further research. Much of the investigations carried out in \cite{OPW,OPW-17a,OPW-17b,PW-2017,PW-2017b,VY2013}
on ${\mathcal U}(\lambda)$ and some other related classes could be considered for further research with an aim to obtain meromorphic
analogue of these classes.
%
}
\er

\subsection*{Acknowledgments}
The authors thank the referee for useful comments. The first author acknowledged the support from a USM research university grant 1001.PMATHS.8011038.
The  work of the second author is supported by Mathematical Research Impact Centric Support of
Department of Science and Technology (DST), India  (MTR/2017/000367).

%

\end{document}